\documentclass[letter,12pt]{amsart}
\usepackage{subfigure}
\usepackage{amsfonts}
\usepackage{newlfont}
\usepackage[centertags]{amsmath}
\usepackage{amsfonts}
\usepackage{amssymb}
\usepackage{graphicx}
\usepackage{plain}
\usepackage{mathrsfs}
\usepackage{makecell}
\usepackage{enumerate}
\usepackage{tikz}
\usepackage{tikz-3dplot}
\usepackage{float}
\usetikzlibrary{calc}
\usepackage{verbatim}
\tikzset{
	MyPersp/.style={scale=1.8,x={(-0.8cm,-0.4cm)},y={(0.8cm,-0.4cm)},
    z={(0cm,1cm)}},
	MyPoints/.style={fill=white,draw=black,thick}
}

        \newtheorem{theorem}{Theorem}
        \newtheorem{definition}{Definition}
 
 \newtheorem{proposition}{Proposition}
        \newtheorem{lemma}[theorem]{Lemma}
        \newtheorem{corollary}[theorem]{Corollary}
        \newtheorem{remark}[theorem]{Remark}
\newtheorem{example}{Example}[section]        
\usepackage{rotating}

\usepackage[colorinlistoftodos]{todonotes}

\newcommand{\sn}{ |\kern-0.25ex|\kern-0.25ex| }
\usepackage{epsfig}

\newcommand{\RR} {\mathbb{R}}
\newcommand{\HH} {\mathbb{H}}
\newcommand{\CC} {\mathbb{C}}

\title[On Polar actions invariant solutions]{On polar actions invariant solutions of semi-linear equations on manifolds}
\author{E. Becerra \and J. Galvis \and N. Martinez}
\address{Universidad Nacional de Colombia, Departamento de 
Matem\'aticas}
\curraddr{Carrera 30 Calle 45, Ciudad Universitaria, Bogot\'a , Colombia}
\email{esbecerrar@unal.edu.co, nmartineza@unal.edu.co, jcgalvisa@unal.edu.co}
\keywords{Semi-linear PDEs, Polar actions, exponential coordinates}

\begin{document}

\begin{abstract}
In this paper we put together  some tools from differential topology and analysis to study second order semi-linear partial differential equations on a Riemannian manifold $M$.  We look for solutions that are constants along orbits of a given group action.  Using some results obtained by Helgason in [J DIFFER GEOM,6(3), 411-419] we are able to write a (reduced) second order semi-linear problem on a submanifold  $\Sigma$. This submanifold is, in certain sense, transversal to the orbits of the group actions and its existence is assumed.  We  describe precise  conditions on the Riemannian Manifold $M$  and the submanifold $\Sigma$ in order to be able to write the reduced equation on $\Sigma$.  These conditions are satisfied by several particular cases including some examples treated separately in the literature such as the sphere, surfaces of revolution and others.  Our framework also includes the setup of polar actions or exponential coordinates.  Using this procedure,  we are left with a second order semi-linear equation posed on a 
submanifold. In particular, if the submanifold  $\Sigma$ is one-dimensional, we can use suitable tools from analysis to obtain existence and properties of solutions.   
\end{abstract}
\maketitle

\tableofcontents
\section{Introduction}

In this paper we study existence of non-negative solutions of  second order semi-linear problem posed on a $n$-dimensional Riemannian manifold $M$. In particular we consider the problem,
\begin{equation}\label{eq:sl-M}
\Delta_M u(x) +b(x)f(u(x))=0,\quad x\in M,\\ 
\end{equation}
where $\Delta_M$ is the Laplace-Beltrami operator defined on $M$ and $u$ is a non-negative function.  The function $f:\mathbb{R}\to \mathbb{R}$ represents the nonlinearity and $b:M \to \mathbb{R}$ is a regularizing coefficient.  The semi-linear problem \eqref{eq:sl-M} has been studied in several configuration manifolds like $\mathbb{R}^n, S^n$ or some surfaces of revolution; see \cite{XueChao,castro1,brezis1983positive,MR2862087,fischer2014infinitely} among others. On these cases, the existence of solutions have been obtained with an analytical approach after a change of variables. For instance, in the case 
of $\mathbb{R}^n$ the problem is 
reduced to a one-dimensional problem when considering polar coordinates and seeking form radial solutions. The resulting problems is a one dimensional semi-linear second order problem having possible some finite number of singularities. Then, existence of solutions of the resulting one-dimensional problem is obtained using tools from analysis. Solution of the one-dimensional reduced problem corresponds to  solutions of the original semi-linear $n$-dimensional problem that solely depend on the radial part when considering them in polar coordinates. 

Note that the construction  mentioned above, depends on some particular properties the respective configuration space that requires a particular coordinated system. Indeed, what is behind this construction is a particular geometry on the configuration space; in other words, the solution to the problem  depends on the geometry of the configuration space $M$. These examples are particular examples of a family of geometries described via a general polar coordinates system coming from the set of symmetries of the manifold $M$, sometimes referred as {\tt polar actions}. See \cite{Jo,DC}. 

The main goal  of this manuscript is to present a geometrical point of view of some reported solutions methods of the problem \eqref{eq:sl-M}. This geometric interpretation allows us to extended 
the analysis to second order semi-linear problems posed on more general configuration spaces. This approach lies on a characterization of the Laplace-Beltrami operator on manifolds with polar actions; see \cite{HelgasonIII}. More precisely, 
If we consider a polar action of a topological group $G$ on the manifold $M$ with one-dimensional transverse submanifold $\Sigma$, under usual assumptions on regularity and growth
condition on the non-linear term, we obtain that the problem \eqref{eq:sl-M} has at least one non-negative solution $u\in C^\infty(M)$. The 
assumptions on the nonlinear term are classical assumptions; for instance, 
they are similar to a result given in \cite{XueChao}.

We stress that the polar action is the geometric properties of the configuration space we mentioned above. For this action we ask an additional condition on the rank of the transverse submanifold $\Sigma$. This  is a technical condition that appears in order to reduce the problem \eqref{eq:sl-M} posed 
on $M$ to a problem posed on $\Sigma$ when seeking for solutions that are invariant with respect to the polar action of $G$. In this paper we consider mainly the case where $\Sigma$ is a one-dimensional manifold (that is, and interval) since for these case 
there are available results to obtain existence and properties of solutions; 
see \cite{XueChao} and references therein. In the case of $\Sigma$ being of dimension  greater than one, we can still write a second oder semi-linear problem on $\Sigma$ but, of course, obtaining existence and properties of solutions is more involved in this case and it will require more sophisticated tools that 
the ones used to obtain existence and properties 
of solutions in the one-dimensional case. The case where $\Sigma$ is of dimension larger that one, will be object of future research. 

The organization of the paper is the following: in Section~\ref{sec:special solutions} we give a general method to solve the problem in the cases of the euclidean space, the sphere and surface of revolution. Indeed, the semi-linear problem on each case can be reduced to the same type of one-dimensional problem (compare the ordinary differential equations \eqref{eq:equivPDE-Rn}, \eqref{eq:equivPDE-Sn} and \eqref{eq:equivPDE-SR}). In Section~\ref{sec:geometry} we present a brief summary of the geometric tools used to reduce the problem \eqref{eq:sl-M} posed on $M$ to a 
second order semi-linear problem posed on the submanifold $\Sigma$.  In particular we give a description of the Laplace-Beltrami operator and we introduce the polar actions. Also, the case of the two-point homogeneous spaces is described with special interest. Finally, Section~\ref{sec:solution} is devoted to give the proof of the main result, Theorem \ref{thm:main}, which gives existence of solution to \eqref{eq:sl-M} via general polar coordinates.

\section{Preliminary examples}\label{sec:special solutions}
Here we  present a standard method to obtain solutions of problem \eqref{eq:sl-M} for the particular cases of the euclidean space, the sphere and surfaces of revolution. The method we present here can be summarized as follows:
\begin{enumerate}
\item Once we describe the Laplace-Beltrami operator we use a change of variables (to some suitable coordinates) to get an equivalent lower dimensional problem.
\item We restrict us to seek for an ansatz  that depends only on one parameter of the change of variable which allows as to reduce the dimension of  the domain where the problem is posed. Here we have to assume that all the coefficients and 
of the differential equation depend solely 
on this one dimensional parameter.
\item by using theory  of semi-linear equations (see for example Theorem~1.1 of \cite{XueChao}) we get conditions that guarantee the solution of the obtained one-dimensional problem. 
\end{enumerate}

\begin{example}[The euclidean space $M=\mathbb{R}^n$]\label{ex:Rn}
Suppose that the functions $u$ and $b$ are radial on polar coordinates $(r,\theta_1,\dots,\theta_{n-1})$. The equation in \eqref{eq:sl-M} becomes into
\[  u'' + \frac{n-1}{r} u'+b(r)f(u)=0,\]
which has a singularity at $r=0$. Under the change of variable  $s=\int_{1}^r c_nt^{1-n} dt$ we get the following equivalent one-dimensional problem
\begin{equation}\label{eq:equivPDE-Rn}
z''(s)+(r(s))^{2(n-1)}b(r(s))f(z)=0
\end{equation}
which is posed on $\mathbb{R}$. $\diamond$
\end{example}
\begin{remark}
Note that the same argument, but for $r$ belonging to the interval $(R_1,R_2)$, can be used to write a similar one-dimensional problem in the case when $M$ is the annulus bounded by the spheres of radius $R_1$ and $R_2$.
\end{remark}
\begin{example}[The  sphere $M=S^{n}$]\label{ex:Sn}

Suppose that the functions $u$ and $b$ are radial on spherical coordinates $(r,\varphi_1,\dots,\varphi_{n-1})$ where $r$ is the arc-length of a meridian from a fixed point $p\in S^n$. The equation \eqref{eq:sl-M} becomes into
\[u'' + (n-1)\frac{\cos(r)}{\sin(r)} u'+b(r)f(u)=0\]
which has a singularity at $r=0$ and 
$r=\pi$. Under a suitable change of variable and defining $s=\int_{r_0}^r c_n\sin(t)^{1-n} dt$ we get the following equivalent ODE
\begin{equation}\label{eq:equivPDE-Sn}
z''(s)+(c_n\sin(r(s))^{n-1})^2b(r(s))f(z)=0.
\end{equation}
 
\begin{remark}
If we consider the case $n=2$, $S^2\subset \mathbb{R}^3$ and $p=(0,0,1)$ and the north pole. We can 
change variable to the vertical axis $z$ using $\sin(r)=\sqrt{1-z^2}$ to get the equation
\[
(1-z^2)  \frac{d ^2v}{dz^2} + (1-n)z\frac{dv}{dz}+b(z)f(v)=0.
\]
where $v(z)=u(r)$.
Note that due to the change of variables now the singularities are  located at $z=-1$  and $z=1$.
This case was studied in \cite{fischer2014infinitely} where they show that the semi-linear Laplace-Beltrami equation has infinitely many solutions on the unit sphere which are symmetric with respect to rotations around some axis. They use a fixed point theorem combined with some energy analysis to obtain local solutions. To obtain global solutions they write a Pohozaev-type identity to prove that the energy is continuously increasing with the initial condition and then use phase plane analysis to prove the existence of infinitely many solutions.
\end{remark} $\diamond$
\end{example}

Under positivity and regularity assumptions on $b$ and asymptotic behavior of $f$, the equations~\eqref{eq:equivPDE-Rn} has regular and positive solution $z(s)$ whenever positive solutions of $z''+b(r(s))(r(s))^{2(n-1)}=0$ exists. Same result holds if we consider the equation~\eqref{eq:equivPDE-Sn} and positive solutions of $z''+b(r(s))(c_n\sin(r(s))^{n-1})^2=0$ also exists.

On both cases, we get conditions on $b$ and $f$ in order to guarantee solution to the problem \eqref{eq:sl-M}.

\begin{example}[Closed surfaces of revolution]\label{ex:surf.rev.}
Consider now a plane curve $\gamma(t)=(x(t),z(t))$ defined on an interval $I=[a,b]$ with non-negative components and parametrized by arc-length. If in addition we assume that $x(a)=x(b)=0,$ then we can construct a simply connected surface of revolution $S$ parametrized by $(x(t)\cos(\theta),x(t)\sin(\theta),z(t))$. The surface $S$ has a Riemannian structure coming from the euclidean structure on $\RR^3$ and it satisfies that the meridian geodesics with origin on $p=(0,0,z(a))$ are the curves with $\theta$ constant \cite{DC}. Indeed, the coordinates $(t,\theta)$ are the geodesic polar coordinates and the Laplace-Beltrami operator on geodesics radial functions, using \eqref{eq:LB-polar.coord.}, is 
$$\Delta u=u''(t)+(\ln(x(t)))'u'(t),$$
and the problem \eqref{eq:sl-M} is written as
\begin{equation}
(x u')'+b\varphi f (u)=0.
\end{equation}
Again, by using a change of variable of the type $s=\int_{r_0}^r \frac{1}{A(t)} dt$ where $A(t)$ is the area of the geodesic sphere of radius $t$, we obtain that the previous system is given by 

\begin{equation}\label{eq:equivPDE-SR}
z''(s)+(A(r(s))^2b(r(s))f(z)=0.
\end{equation}
which is an ODE.$\diamond$
\end{example}

The remarkable fact is that the solutions of the problem \eqref{eq:sl-M} on these three cases can be studied by the same type of second order ODE (cf. equations \eqref{eq:equivPDE-Rn},\eqref{eq:equivPDE-Sn} and \eqref{eq:equivPDE-SR}). Same procedure can be applied to more general situations where we can find a ''radial" variable and the remaining variables can be considered as ''symmetries"-variables.

\section{Geometrical setting}\label{sec:geometry}
This section is devoted to present some basic facts on Riemannian geometry needed to describe our main problem in a general set-up. The key ingredients will be the Laplace-Beltrami operator and the polar actions.

\subsection{Riemannian geometry}
 Let us begin by presenting some familiar concepts. For a more detailed description we reefer the reader to~\cite{Jo}. 
 
A manifold $M$ of dimension $n$ is a connected para-compact Hausdorff space for which every point $x\in M$ has a neighborhood $U_x$ that is homeomorphic to an open subset $\Omega_x$ of $\RR^n$. Such a homeomorphism $\phi_x : U_x \to \Omega_x$ is called a {\it coordinate chart}. If for two charts the function $\phi_x\circ \phi_y^{-1}:\Omega_y\to \Omega_x$ is a $C^r$-diffeomorphism we say that the manifold has a $C^r$-differentiable structure.  We denote $T_xM$, the vector space which consist of all tangent vector to curves in $M$ on the point $x$. It is called the \emph{tangent space} of $M$ at the point $x$. 

A {\it Riemannian metric} on a differentiable manifold $M$ is given by a scalar product on each tangent space $T_x M$ which depends smoothly on the base point $x$. A {\it Riemannian manifold} is a differentiable manifold equipped with a Riemannian metric.

In any system of local coordinates $(x_1,\ldots,x_n)$ the Riemannian metric is represented by the positive definite, symmetric matrix $(g_{ij}(x))_{i,j=1,\ldots,n}$ where the coefficients depend smoothly on $x$.  

Let $\gamma : [a, b]\to  M$ be a smooth curve. The length of $\gamma$ is defined as $$L(\gamma)=\int_a^b \|\gamma'(t)\|dt$$ where the norm of the tangent  vector $\gamma'$ is given by the Riemannian metric. This value is invariant under re-parametrization. Taking the infimum of the values $L(\gamma)$ among all the curves $\gamma$ joining two points $p,q\in M$ we can define a distance function on $M$ and the topology of this distance coincides with the topology of the manifold structure of $M$.
\begin{theorem}
Let $M$ be a Riemannian manifold, $x\in M$ and $v\in T_xM$. Then there exist $\epsilon>0$ and precisely one geodesic $c : [0, \epsilon] \to M$ with $c(0) = x, c(0)= v$. In addition, $c$ depends smoothly on $x$ and $v$.
\end{theorem}
The main consequence of this theorem is that we can define what is called {\it normal coordinates} via exponential map. The {\it exponential map} is a diffeomorphism between an open set of $T_xM$ (with center at $0$) defined by geodesics and an open on $M$ (with center at $x$).  When we use standard coordinates $(r, \varphi)$, where $\varphi = (\varphi_1, \ldots,\varphi_{n-1})$ parametrizes the unit sphere $S^{n-1}$ on $\RR^n$, we then obtain polar
coordinates on $T_xM$ (via an orthonormal linear isomorphism $T_xM\equiv \RR^n$), thus we get new coordinate system on $M$ via the exponential map. Such coordinates are known as {\it geodesic polar coordinates} and satisfies that the lines with $\varphi$ constant are geodesic.

\begin{corollary}
For any $x\in M$, there exists $\rho>0$ such that geodesic polar coordinates can be introduced on $B(x, \rho) := \{q \in M : d(x, q) \leq \rho\}$. For any such $\rho$ and any $q \in \partial B(x, \rho)$, there is precisely one geodesic of shortest length $(= \rho)$ from $x$ to $q$, and in polar coordinates, this geodesic is given by the straight line $x(t) = (t, \varphi_0 ), 0 \leq t \leq \rho$, where $q$ is represented by the coordinates $(\rho, \varphi_0), \varphi_0 \in S^{n-1}$.
\end{corollary}

%
%
\begin{example}[The sphere]
The surface $S^2$ can be constructed via radial circles $C_z$ centered at some point in $[-r_0,r_0]$ in the $z$-axis. It provides us a parametrization of $S^2$ via the formula  $$\gamma(z,t)=(\sqrt{r_0^2-z^2}\cos(t), \sqrt{r_0^2-z^2}\sin(t),z)$$ for a fixed value of $r_0$. 
We consider the tangent space $T_pS^2$ at some point $p\in S^2$ (see FIGURE 1). In $T_pS^2$ we can define the (planar) polar  coordinates $(r,\theta)$  centered at the point $p$. Without loss of generality we assume $p=(0,0,r_0)$ and we can project each circle defined on these polar coordinates to $S^2$ such that the image of these circles coincides with the circles $C_z$ and the radio $r$ projects on a geodesic transversal to the circles $C_z$.
\tdplotsetmaincoords{80}{125}
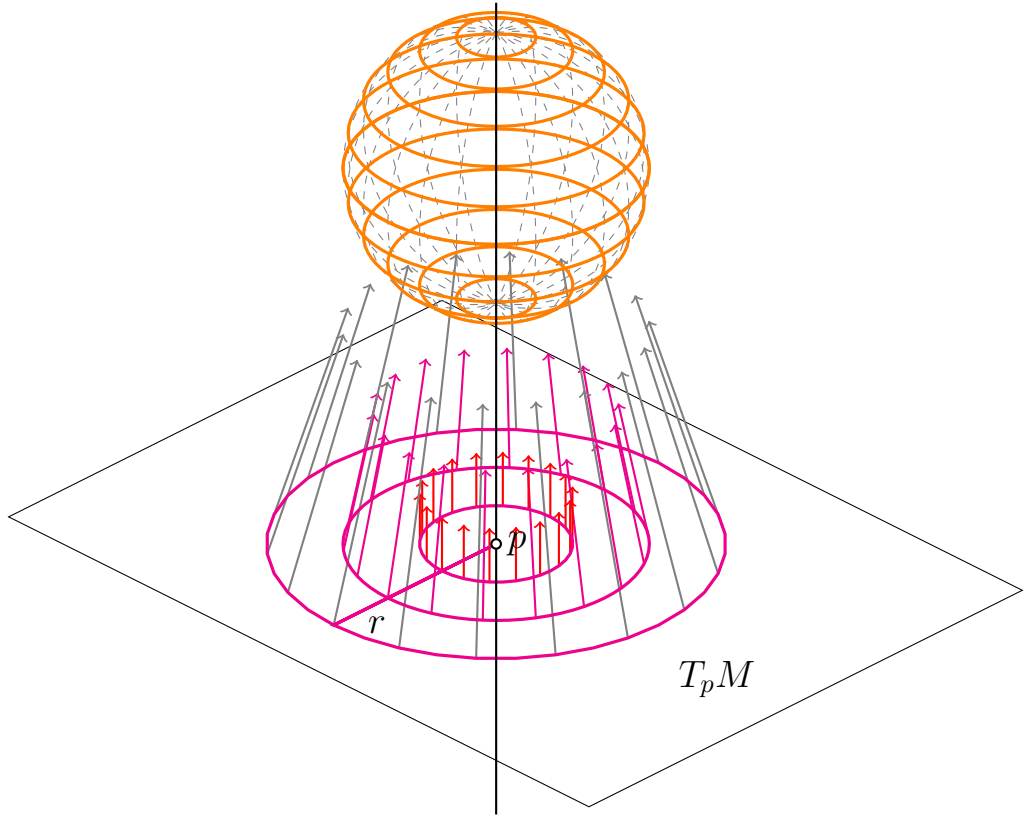
\begin{figure}[t]
\centering
\begin{tikzpicture}[MyPersp,font=\large]
	\def\h{4}
	\def\a{35}
	\def\aa{35}
	\pgfmathparse{\h/tan(\a)}
  \let\b\pgfmathresult
	\pgfmathparse{sqrt(1/cos(\a)/cos(\a)-1)}
  \let\c\pgfmathresult 
	\pgfmathparse{\c/sin(\a)}
  \let\p\pgfmathresult 
	\coordinate (A) at (2,0.5*\b,0);
	\coordinate (B) at (-2,0.5*\b,0);
	\coordinate (C) at (-2,-1.5,{(1.5+\b)*tan(\a)});
	\coordinate (D) at (2,-1.5,{(1.5+\b)*tan(\a)});
	\coordinate (E) at (2,-2.5,0);
	\coordinate (F) at (-2,-2.5,0);
	\coordinate (CLS) at (0,0,{\h-\p});
	\coordinate (CUS) at (0,0,{\h+\p});
	\coordinate (FA) at (0,{\c*cos(\a)},{-\c*sin(\a)+\h});
	\coordinate (FB) at (0,{-\c*cos(\a)},{\c*sin(\a)+\h});
	\coordinate (SA) at (0,1,{-tan(\a)+\h}); 
	\coordinate (SB) at (0,-1,{tan(\a)+\h});
	\coordinate (PA) at ({sin(\aa},{cos(\aa)},{\h+3*\p});
	\coordinate (PB) at ({sin(\aa},{cos(\aa)},{\h-3*\p});
	\coordinate (P) at ({sin(\aa)},{cos(\aa)},{-tan(\a)*cos(\aa)+\h});

	\draw (E)--(A)--(B)--(F)--cycle;


	\foreach \t in {20,40,...,360}
		\draw[->,magenta,thick] ({cos(\t)},{sin(\t)},0)
     --({0.8*cos(\t)},{0.8*sin(\t)},{0.25*\h});
     \foreach \t in {20,40,...,360}
		\draw[->,red,thick] ({0.5*cos(\t)},{0.5*sin(\t)},0)
     --({0.5*cos(\t)},{0.5*sin(\t)},{0.1*\h});
     \foreach \t in {20,40,...,360}
		\draw[->,gray,thick] ({1.5*cos(\t)},{1.5*sin(\t)},0)
     --({cos(\t)},{sin(\t)},{0.4*\h});
	\draw[magenta,very thick] (1,0,0) 
		\foreach \t in {5,10,...,360}
			{--({cos(\t)},{sin(\t)},0)}--cycle;
            \draw[magenta,very thick] (1,0,0) 
		\foreach \t in {0,10,...,360}
			{--({0.5*cos(\t)},{0.5*sin(\t)},0)}--cycle;
            \draw[magenta,very thick] (1,0,0) 
		\foreach \t in {0,10,...,360}
			{--({0.01*cos(\t)},{0.01*sin(\t)},0)}--cycle;
            \draw[magenta,very thick] (1,0,0) 
		\foreach \t in {0,10,...,360}
			{--({1.5*cos(\t)},{1.5*sin(\t)},0)}--cycle;

	\foreach \i in {-1}{
		\foreach \t in {0,15,...,165}
			{\draw[gray,dashed] ({cos(\t)},{sin(\t)},\h+\i*\p)
				\foreach \rho in {5,10,...,360}
					{--({cos(\t)*cos(\rho)},{sin(\t)*cos(\rho)},
          {sin(\rho)+\h+\i*\p})}--cycle;
			}
		\foreach \t in {-75,-60,...,75}
			{\draw[orange,very thick] ({cos(\t)},0,{sin(\t)+\h+\i*\p})
				\foreach \rho in {5,10,...,360}
					{--({cos(\t)*cos(\rho)},{cos(\t)*sin(\rho)},
         {sin(\t)+\h+\i*\p})}--cycle;
			}
					\draw[orange,very thick] (1,0,{\h+\i*\p})
		\foreach \t in {5,10,...,360}
			{--({cos(\t)},{sin(\t)},{\h+\i*\p})}--cycle;
		}
	\draw[black,thick] (0,0,-2)--(0,0,4);
	\fill[MyPoints] (0,0,0) circle (1pt)node[right]{$p$};
	\fill[MyPoints] (0,2.5,0) node[left]{$T_pM$};
	\fill[MyPoints] (1.3,0.2,0) node{$r$};
\end{tikzpicture}
\label{figura1}\caption{Polar coordinates on $S^2$}
\end{figure}
\end{example}

\begin{example}[Two points homogeneous spaces]\label{two points homogeneous spaces}
A Riemannian manifold $(M,g)$ is {\it two points homogeneous spaces} if the group $I(M)$ of isometries acts transitively on the space of equidistance pair of points. The definition means that for any $p_1,p_2,p_1',p_2'\in M$ so that $d_g(p_1,p_2)=d_g(p'_1,p'_2)$ then there exists $\varphi\in I(M)$ so that $\varphi(p_i)=p_i'$.  As a direct consequence of the definition we have that the function  $A_p(r)$ defined as the Riemann measure of the geodesic sphere $S_r(P)$ with center at $p\in M$ is independent of the point $p$ hence $A$ is globally defined on $M$.

Note that from definition we can observe that we only need one invariant to determine all the equidistant pair of points. The last statement is just the well known fact that any two points homogeneous spaces is the euclidean space or a symmetric space of rank one (see Chapter X in \cite{Hel}). For the complete classification of the symmetric spaces of arbitrary rank we refer to \cite[Table V Chapter X]{Hel2}, and for the convenience of the reader we give the list of the two points homogeneous spaces. 
\begin{itemize}
\item The euclidean space $\RR^n$ and the spheres $S^n$ for $n\geq 1$
\item Real, complex and quaternion projective spaces $P\RR^n, P\CC^n,$ and $P\HH^n $ for $n\geq 2$
\item Real, complex and quaternion hyperbolic space $H\RR^n, H\CC^n,$ and $H\HH^n $ for $n\geq 2$
\item Cayley projective $P\CC_a^2$ and hyperbolic spaces $H\CC_a^2$.
\end{itemize}
Some special cases of closed surfaces of revolution, classified by the Gaussian curvature, are isometric with two points homogeneous spaces.
\end{example}

\subsection{The Laplace-Beltrami operator and geodesic polar coordinates}
What we need now is the description of the Laplace-Beltrami operator for Riemann manifolds $M$. If we denote by $\bar{g}=\det(g_{ij})$, the Laplace-Beltrami operator acting on a smooth function
$u:M\to \mathbb{R}$ is defined as\footnote{\mbox{ For the coordinate-free expression of $\Delta_Mu$ we refer to \cite[Sec~2.1]{Jo}.}} 
$$\Delta_Mu=\frac{1}{\sqrt{\bar{g}}}\sum_{k}\frac{\partial}{\partial x_k}\left( \sum_{i} g^{ik}\sqrt{\bar{g}}\frac{\partial}{\partial x_i}u \right)$$ where $(g^{ij}):=(g_{ij})^{-1}$.

In the case of geodesic polar coordinates we have a special form of the previous operator:
\begin{equation}\label{eq:LB-polar.coord.}
\Delta_M u=\dfrac{\partial^2u}{\partial r^2}+ \dfrac{1}{\sqrt{\bar{g}}}\dfrac{\partial \sqrt{\bar{g}}}{\partial r} \dfrac{\partial u}{\partial r}+ \dfrac{1}{\sqrt{\bar{g}}}\sum_{i,j=1}^{n-1}\dfrac{\partial }{\partial \varphi_i}(g^{ij}\sqrt{\bar{g}}\dfrac{\partial u}{\partial \varphi_j}).
\end{equation}

For two points homogeneous spaces the previous formula takes the following form,
\begin{equation}\label{eq:LB-2pths}
\Delta u=u''+(\ln A)'u'
\end{equation}
where $u$ is a geodesical radial function and $A(r)$ is the Riemann measure of the geodesic sphere of radius $r$ (cf. Lemma 7.12 in \cite{Hel}).

\subsection{Generalized polar coordinates}
In order to give a generalization of the geodesic polar coordinates given in the previous section for two points homogeneous spaces to a more general set-up, we shall do some remarks on the presentation of the examples \ref{ex:Rn}-\ref{ex:surf.rev.}. The geometry of these will serve to us as a motivation to study the more general notion, the \emph{polar coordinates}. For that goal some definitions are needed.
\begin{definition}
An \emph{action} of a Lie group $G$ with identity element $e$ on a manifold $M$ is a differentiable map $\phi:G\times M\rightarrow M$ such that
\begin{itemize}
\item $\phi(e,x)=x$ for all $x\in M$,
\item $\phi(\alpha\beta,x)=\phi(\alpha,\phi(\beta,x))$ for all $\alpha,\beta\in G$ and for all $x\in M$.
\end{itemize}
\end{definition}
For short we denote for all $\alpha\in G$ and $x\in M$: $\phi(\alpha,x)=\alpha x$ and we denote by $$G\cdot x:=\{\alpha x|\alpha\in G\}$$ the orbit of $x\in M$. For a Riemannian manifold $(M,g)$ we say that the action of $G$ on $M$ is an \emph{isometry} if for any $\alpha\in G$ we have $g_p(X,Y)=g_{\alpha p}(d_p\phi(\alpha) X,d_p\phi(\alpha) Y)$ for all $X,Y\in T_pM$. 

With this definition we proceed to revisit some examples. Let us begin with the example of the 2-sphere. The surface $S^2$ can be constructed via radial circles $C_z$ centered at some point in $[-r_0,r_0]$ in the $z$-axis. It provides us a parametrization of $S^2$ via the formula  $$\gamma(z,t)=(\sqrt{r_0^2-z^2}\cos(t), \sqrt{r_0^2-z^2}\sin(t),z)$$ for a fixed value of $r_0$.

We define an action from the compact group $S^1$ on $S^2$ by $$\lambda\cdot \gamma(z,t)=(\sqrt{r_0^2-z^2}\cos(t+\theta_{\lambda}), \sqrt{r_0^2-z^2}\sin(t+\theta_{\lambda}),z),$$ where $\theta_{\lambda}$ is the angle corresponding to $\lambda\in S^1$. Note that the orbit of the point $\gamma(z,t)$ is exactly the circle $C_z$. Now we fix a point $p\in S^2$ (that is not the north neither the south pole) and consider the meridian $\Sigma$ through $p$, we obtain that  $\Sigma=\{(\sqrt{r_0^2-x^2},0,z)|-r_0< z< r_0\}$ for a suitable $x\in (-r_0,r_0)$. For each $z\in (-r_0,r_0)$, the tangent space at $p$ can be decomposed as $T_pS^2=E_1\oplus E_2$ where $E_1$ is a one-dimensional vector space generated by a vector in the direction of $\Sigma$ and $E_2$ is generated by a vector in the direction of the circle $C_z$. Note that in this case the submanifold $\Sigma$ is \emph{transversal} to the action of $S^1$ on $S^2$.  We stress that this fact only holds for points that are not the north or the south pole (i.e. $z\in (-r_0,r_0)$) because this particular points are fixed point of the action; nevertheless such points can be considered in the situation described before by using another (invariant by rotation) parametrization.

This situation also happens in the other examples presented in section \ref{sec:special solutions}. That is, there exists an action by isometries of a compact group $G$ on the Riemannian manifold $M$, where the orbit of a single point is transversal to a fixed submanifold $\Sigma$ of $M$:
\begin{itemize}
\item For  example \ref{ex:Rn}, there is an action of the group $SO(n)$. For this action the submanifold $\Sigma$ can be chosen as an infinite line from the origin in $\RR^n$
\item  For the example \ref{ex:Sn}, the subgroup $G$ of $SO(n+1)$ defined by the rotations about the $x_{n+1}$-axis acts on $S^{n}$ by considering it as a submanifold of the space $\RR^{n+1}$ (that is, $G=SO(n-1)$). The submanifold $\Sigma$ can be chosen as a geodesic line joining the points $(0,\dots ,0, 1)$ and $(0,\dots ,0, -1)$ of $S^n$.
\item For the surfaces of revolution, the group $S^1$ acts by rotations on the surface, where $\Sigma$ is image of the curve $\gamma(t)=(x(t),0,z(t))$. Moreover note that this holds for any surface of revolution not only the closed ones.
\item In general, for two points homogeneous spaces it always exists such transversal one-dimensional submanifold
\end{itemize}

All these considerations also hold for two points homogeneous spaces. In this way, the polar geodesics coordinates on these spaces can be understood as the one given by the orbits of an action and the transversal one-dimensional submanifold.
to a general set-up of actions from a group of isometries $G$ on a Riemannian manifold $M$, where the submanifold $\Sigma$ can be chosen to be \emph{transversal} to the orbits of the action. 

\begin{definition}\label{def:polar action}\label{polaraction}
Let $G$ be a Lie group acting on $M$ by isometries. One says that the action is \emph{polar} if there is a complete immersed submanifold $\Sigma$ in $M$ that is transversal to non trivial orbits. That is, if $G\cdot x\neq\{x\}$ then $\Sigma$ intersects $G\cdot x$ in a single point and
$$T_xM=T_x\Sigma\oplus T_x(G \cdot x).$$\footnote{There is a weaker notion but we do not consider it because we want to avoid fixed points of the $G$-action} 
\end{definition}
From the previous work we can verify that the action of $S^1$ on the sphere $S^n$ is a polar action when consider $\Sigma$ as a meridian without the north and south poles (these are the fixed points of the action). The use of polar action gives meaningful advantages; for instance we can consider more situations than two points homogeneous spaces where we have a suitable polar coordinates, namely, the space $T^2=S^1\times S^1$ with an action of $S^1$ defined by rotating the second coordinate, the set $\{(s_1,s_0): s_1\in S^1\}$ is a submanifold of $T^2$, transversal to the orbits of the action.
The main advantage we work with is that for manifolds which admits polar actions we obtain a general description of the Laplace-Beltrami operator. In particular, we extend the expression of  Laplace-Beltrami operator for general surface of revolution that are not isometric to two points homogeneous spaces.    

\subsubsection{The radial part of a differential operator}
For this section we will assume a polar action and our main reference is \cite{HelgasonIII}. We remake some constructions presented there in order to obtain a clear presentation. 

Let $(M,g)$ be a complete Riemannian manifold and let $\Sigma$ a submanifold of $M$. For any element $s\in \Sigma$, we consider $\Sigma_s^{\bot}$ the set of all the geodesics starting at $s$ transversally to $\Sigma$. For a fixed neighborhood $U_0$ of a point $s_0\in\Sigma$ we can assume that all the $\Sigma_s^{\bot}$ for $s\in U_0$  are disjoints, then, their union defines a neighborhood $V_0$ of $s_0$ in $\Sigma$. For any $C^{\infty}$-function $u$ with compact support defined on $\Sigma$, we define $\tilde{u}$ on $V_0$ as $u(s)$ in $\Sigma_s^{\bot}$, for all $s\in U_0$. Let $D$ a differential operator defined on $M$. We define $D'$ the \emph{projection} of $D$ on $\Sigma$ by the equation $$D'(u)(s)=D(\tilde{u})(s).$$
\begin{proposition}[\cite{HelgasonIII}, Proposition 1.1]
Let $L_M$ and $L_{\Sigma}$ be the Laplace-Beltrami operators on $M$ and $\Sigma$ respectively. Then $L_{\Sigma}$ is the projection of $L_M$.
\end{proposition}
Also we have a local description that allows to define the transversal part of a differential operator;
\begin{lemma}[\cite{HelgasonIII}, Lem1.2] 
Let $\Sigma$ be a transversal submanifold in $M$. For each $x\in M$ there exist a cross section $B\times W$, with $B$ a relatively compact submanifold of $G$ and $W$ and open neighborhood of $x$ in $\Sigma$, such that the mapping $\eta:B\times W\rightarrow V$ defined by $(b,w)\mapsto bw$ is a diffeomorphism onto an open neighborhood $V$ of $x$ in $M$.
\end{lemma}
For $x\in M$, let us consider $N=G\cdot x$ and we construct the neighborhood $N_x^{\bot}$ as before. We can assume $N_x^{\bot}$ is transversal in $M$ and we apply the lemma above to construct the neighborhood $V=\cup_{w\in W} Bw$ of $x$. For $u$ a $C^{\infty}$-function on $M$, we restrict it to $W\subset N_x^{\bot}$ and then we extend this restriction to a function on $V$ defined as $u_x(bw)=u(w)$ for $b\in B$ and $w\in W$. For a differential operator $D$ on $M$ we define the \emph{transversal part} of $D$ by:
$$D_T(u)(x)=D(u_x)(x).$$We denote it by $D_T$. 
\begin{lemma}[\cite{HelgasonIII}, Theorem 1.4 ]
Let $M$ a Riemannian manifold and $G$ a Lie group acting by isometries on $M$. Let $N$ be any orbit by the action of $G$ in $M$ and let $u^-$ denote the restriction to $\Sigma$ of a $C^{\infty}$-function $u$ on (a subset of) $M$. Then$$L_M(u)=L_\Sigma(u^-)+L_T(u)^-$$
\end{lemma}
\begin{theorem}[\cite{HelgasonIII}, Proposition 2.1]\label{teoradialpart}
Let us suppose $\Sigma$ transversal to the action of $G$. Let $D$ be a differential operator on $M$. There is a unique operator $\Delta(D)$ on $\Sigma$ such that 
$D(u)^-=\Delta(D)(u^-)$, for each $G$-invariant function $u$ defined on an open subset of $M$. 
\end{theorem}
The operator $\Delta(D)$ is called the \emph{radial part} of $D$. \\
Henceforth, we suppose $G$ a compact  group of isometries of $M$.
\begin{theorem}[\cite{HelgasonIII}, Theorem 2.11] For any one-dimensional submanifold $\Sigma\subset M$, transversal to the action of $G$, we have
$$\Delta(L_M)=\frac{1}{\sqrt{A(r)}}L_{\Sigma}(\circ\sqrt{A(r)})-\frac{1}{\sqrt{A(r)}}L_{\Sigma}(\sqrt{A(r)})$$ where $A(r)$ denotes the 
Riemannian $G$-invariant measure of the geodesic sphere of radius $r$ along
the submanifold $\Sigma$ and $\circ$ denotes the composition of the operator $L_{\Sigma}$ and the operator multiplication by $\sqrt{A(r)}$. 
\end{theorem}
Let $u$ be a function such that $u(x)=u(gx)$ for all $g\in G$ and for all $x\in M$. Obviously the function $u$ is $G$-invariant, then we can apply the Theorem \ref{teoradialpart} and obtain 
\begin{equation}\label{eq:reduced LB}
L_M(u)=\Delta(L_M)(u)=\frac{1}{\sqrt{A(r)}}L_{\Sigma}(\sqrt{A(r)}u)-\frac{1}{\sqrt{A(r)}}L_{\Sigma}(\sqrt{A(r)})u,
\end{equation}
the explicit formula for the Laplace-Beltrami operator applied on $u$.\\
For some straightforward calculations, the equation \ref{eq:reduced LB} turns out to be\\
\begin{align*}L_M(u)&=\frac{2\nabla(\sqrt{A(r)})\cdot\nabla(u)}{\sqrt{A(r)}}+L_{\Sigma}(u)\\
&=\frac{(A(r))'u'}{A(r)}+u''\\&=(\ln(A))'u'+u'',\end{align*}\\
recovering the equation \eqref{eq:LB-2pths} presented before for two points homogeneous spaces (see Section \ref{two points homogeneous spaces}) in a general setting\footnote{We have to remark that we are actually using the fact that  $\Sigma$ is a one-dimensional submanifold of $M$. }.\\ 

We finish this section by discussing the case of higher dimensional transverse manifold. Note that the Laplace-Beltrami operator in ~\eqref{eq:LB-2pths} , for two points homogeneous spaces, applied to the PDE in the problem \eqref{eq:sl-M} yields the following equivalent problem for radial functions 
\begin{equation}
(Au')'+Abf(u)=0.
\end{equation}
This one dimensional equation correspond to a general situation presented in the examples leading to equations \eqref{eq:equivPDE-Rn},\eqref{eq:equivPDE-Sn} and \eqref{eq:equivPDE-SR}. 

In the general case of polar actions, the dimension of $\Sigma$ could be greater than 1, and in these cases the PDE in \eqref{eq:sl-M} is equivalent to 
\begin{equation}\label{eqdimg}
\Delta u+2\dfrac{\nabla \sqrt{A(r)}\cdot \nabla u}{\sqrt{A(r)}}+b(r)f(u)=0
\end{equation}
for $G-$invariant function $u$. Hence we get a semi-linear PDE that can be considered for analysis. In particular, several questions could, in principle be studied:
\begin{itemize}
\item Existence of solutions (that will imply  existence of polar action invariant solution of the original problem)
\item Existence of infinitely many polar actions invariant solutions since the reduced dimensional problem may be parametrized by some boundary data
\item Existence of solution with some desirable properties such as non-negativity and similar properties
\end{itemize}
and many other analytical questions.

To illustrate the use of the geometrical tool presented so far 
we present, in the next section, the case where the parameter $r$ is one dimensional and we are seeking the existence of a non-negative solution of problem \eqref{eq:sl-M}. The general case will be object of future work by some of the authors.

\section{Non-negative solutions for the case of one dimensional transverse submanifold}\label{sec:solution}
We start by stating our result.
\begin{theorem}\label{thm:main}
Consider a polar action of a topological group $G$ on the manifold $M$ with 1-dimensional transverse submanifold $\Sigma$ and denote by $\phi:\Sigma\to (0,\infty)$ the Riemann measure of the geodesic sphere. Let $r_0\in \Sigma$. In addition assume that
\begin{enumerate}[(i)]
\item the change of variables 
 $s=J(r):=\int_{r_0}^r \phi(t)^{-1}dt$ maps $\Sigma$ to $\RR$
\item the linear problem $$z''(s) +b(r(s))\phi(r(s))^2=0,\quad s\in\RR $$ with 
$b(r(\cdot))$ an $\alpha$-regular positive function, has a unique positive $C_0(\RR)\cap C^{2+\alpha}_{loc}(\RR)$ solution,
\item $f\in C^1((0,\infty), (0,\infty))$ is so that $\lim_{q\to 0^+} f(q)/q=\infty$ and we have $\lim_{q\to \infty} f(q)/q=0$.
\end{enumerate}
Then the problem \eqref{eq:sl-M} has at least one non-negative solution $u\in C^\infty(M)$. 
\end{theorem}

Assumption $(i)$ is a technical assumption and it refers to the 
non-integrability of the function 
$1/\phi$. The case where the change of variables $J$ maps $\Sigma$ to a bounded interval is easier to deal with and it is not presented here. The other to hypothesis (ii)-(iii) recall a result given in \cite{XueChao} and the polar action is the geometric properties of the configuration space we mentioned above.  The other condition, the rank of the transverse submanifold, is a technical condition that appears in order to reduce the problem \eqref{eq:sl-M} posed 
on $M$ to a problem posed on $\sigma$ when seeking for solutions that are invariant with respect to the polar action of $G$. 

\begin{proof}
Here we show the solution for the problem \eqref{eq:sl-M}
where $M$ is a surface of revolution or two points homogeneous spaces. As we 
mentioned in the previous section, if we assume  $u$ and $b$ radial functions (in global polar coordinates), then the problem transforms to
$$(SL1)\begin{cases}
(\phi(r)u'(r))' +b(r)\phi(r)f(u(r))=0,\\
u(0)=d,  u(x)\geq 0, 
\end{cases}$$
where $d$ is a real number. Note that  $\phi(r)=x(t)$ in the surface of revolution case and $\phi(r)=A(r)$ in the  two points homogeneous spaces case. On both cases the function $\phi(r)$ is non-negative for $r\neq 0$, so we can define $s=J(r)=\int_{r_0}^r \phi(t)^{-1}dt$ with $r_0>0$  and $z(s)=u(J^{-1}(s))$. It is routine to verify that
$$\dfrac{d r}{ds}=\phi(r) \hspace{1cm}\mbox{\ and\ }\hspace{1cm} \dfrac{dz}{ds}=u'(r)\phi(r),$$

thus $(SL1)$ turns to be equivalent to 
$$(SL2)\begin{cases}
z''(s) +b(r(s))\phi(r(s))^2f(z(s))=0,\\
 z(r_0)=d,  u(r(s))\geq 0, 
\end{cases}$$
where $d$ is any real number. Following Theorem 1.1 in \cite{XueChao} we can conclude that if $z''(s)+b\phi^2=0$ has positive and $C_0(\RR)\cap C^{2+\alpha}_{loc}(\RR)$ solution with $b$ positive and $\alpha$-regular function, then the problem $(SL2)$ has solution if in addition we suppose that $f\in C^1((0,\infty), (0,\infty))$ is so that $\lim_{q\to 0^+} f(q)/q=\infty$ and $\lim_{q\to \infty} f(q)/q=0$.

Note that in this case $u(r)=z(J(r))$ is solution of the problem \eqref{eq:sl-M}.
\end{proof}

\section{Conclusions}

In this paper we study second oder semi-linear partial differential equations on a Riemannian manifold.  In particular, we prove the existence of  solutions that are constants along orbits of a given group action.  Using some results obtained by Helgason in
\cite{Radialpart} we reduce the dimension of the partial differential equation and are able to bring known tools from analysis to obtain the results.
We detailed all the geometrical constructions needed in order to obtain the reduced dimensional problem that, in the general case, is posed on a submanifold of the original domain where the partial differential equation is posed; see~\eqref{eqdimg}.

\section*{Acknowledgements}
J. G. wants to thank professors A. Castro and G. Rodr\'iguez for useful conversations concerning this type of problems and exponential coordinates.

\end{document}